\def\mapright#1{\smash{\mathop{\longrightarrow}\limits^{#1}}}

\def\mapup#1{\uparrow\rlap{$\vcenter{\hbox{$\scriptstyle#1$}}$}}
\def\mapupleft#1{\uparrow\llap
                   {$\vcenter{\hbox{$\scriptstyle#1$}}\;\;$}}

\documentclass{amsart}
\usepackage{amscd, amsmath, amsthm, amssymb}

\newtheorem{theorem}{Theorem}[section]
\newtheorem{lemma}[theorem]{Lemma}
\newtheorem{proposition}[theorem]{Proposition}

\theoremstyle{definition}     

\theoremstyle{remark}
\newtheorem{remark}[theorem]{Remark}

\numberwithin{equation}{section}

\begin{document}

\title[Alternating group and K3 surfaces]
{Extensions of the alternating group of degree 6 in the geometry
of K3 surfaces}

\author[J.H. Keum]{JongHae Keum}
\address{School of Mathematics, Korea Institute for Advanced Study,
Dongdaemun-gu, Seoul 130-722, Korea}
\email{jhkeum@kias.re.kr}
\thanks {The first author was supported by Korea Research Foundation Grant
(KRF-2002-070-C00003).}

\author[K. Oguiso]{Keiji Oguiso}
\address{Graduate School of Mathematical Sciences,
University of Tokyo, Komaba, Meguro-ku,
Tokyo 153-8914, Japan
}
\email{oguiso@ms.u-tokyo.ac.jp}

\author[D. -Q. Zhang]{De-Qi Zhang}
\address{Department of Mathematics, National University of Singapore,
2 Science Drive 2, Singapore 117543, Singapore}
\email{matzdq@math.nus.edu.sg}


\begin{abstract}
We shall determine the uniquely existing extension of the
alternating group of degree $6$ (being normal in the group) by a
cyclic group of order $4$, which can act on a complex K3 surface.
\end{abstract}

\maketitle


\setcounter{section}{0}
\section{ Introduction}

A K3 surface $X$ is a simply-connected compact complex
$2$-dimensional
manifold which admits a nowhere vanishing holomorphic $2$-form $\omega_{X}$.
As is well known, K3 surfaces form a $20$-dimensional analytic family.

In our previous note [KOZ], we have shown, among all K3 surfaces,
the unique existence of the triplet $(F, \tilde{A}_{6}, \rho_{F})$
of a complex $K3$ surface $F$ and its (faithful) finite group
action $\rho_{F} : \tilde{A}_{6} \times F \longrightarrow F$ of
$\tilde{A}_{6}$ on $F$, up to isomorphisms. The group
$\tilde{A}_{6}$ is an extension of $A_{6}$ (being normal in
$\tilde{A}_{6}$) by $\mu_{4}$, which has also been shown to be the
unique maximum finite extension of $A_{6}$ in the automorphism
groups of $K3$ surfaces. (Here and hereafter, we shall employ the
notation of groups as in the list of notation at the end of
Introduction.) We have also described the target $K3$ surface $F$;
it is isomorphic to the minimal resolution of the surface in
$\mathbf P^{1} \times \mathbf P^{2}$ given by the following
equation, where $([S:T],[X:Y:Z])$ are the homogeneous coordinates
of $\mathbf P^{1} \times \mathbf P^{2}$:
$$S^2(X^3+Y^3+Z^3) - 3(S^2+T^2)XYZ = 0.$$

However, as was remarked in [KOZ], the action of $\tilde{A}_{6}$
is so invisible in the equation that it seems hard
to find the abstract group structure of the (uniquely existing) group
$\tilde{A}_{6}$ from the equation above and it remains unsolved to date.

The aim of this short note is to describe the group structure of
$\tilde{A}_{6}$ explicitly as an abstract group (Theorem (4.1)).

In contrast
to the fact that
$\text{Out}(A_{n}) \simeq C_{2}$ when $n \geq 3$ and $n \not= 6$,
the very special nature
of $A_{6}$, that {\it the outer automorphism
group $\text{Out}(A_{6})$ is isomorphic to a bigger group
$C_{2}^{2}$}, makes our humble work non-trivial and interesting.
Indeed, corresponding to the three involutions,
$\text{Aut}(A_{6})$ has three index $2$ subgroups $A_{6} < G <
\text{Aut}(A_{6})$, which are $S_{6}$, $\text{PGL}(2, 9)$ and
$M_{10}$. (See for instance [Su, Chapter 3], [CS, Chapters 10, 11].
See also Section 1). According to Suzuki [Su, Page 300],
this extraordinary property of $A_{6}$ also makes the classification of
simple groups deep and difficult.

We first notice that there are exactly four isomorphism classes of
$A_{6}.\mu_{4}$ corresponding to the four normal proper subgroups
$A_{6}$, $S_{6}$, $\text{PGL}(2, 9)$ and $M_{10}$ of $\text{Aut}\,
A_{6}$ (Theorem (2.3)). This is purely group-theoretic and should
be known to experts. We then determine which one our
$\tilde{A_{6}}$ is. The part here involves geometric arguments on
K3 surfaces such as the representation of the group action on the
Picard lattice of a K3 surface and its constraint from projective
geometry of a target K3 surface (Proposition (3.2), Lemma (4.2)).
It turns out that our $\tilde{A}_{6}$ is the one which arises from
the last normal subgroup $M_{10}$, the Mathieu group of degree
$10$ (Theorem (4.1)).

It might be of interest that K3 surfaces could distinguish $M_{10}$
from $S_{6}$ and $\text{PGL}(2, 9)$ in this way.
\par \vskip 1pc \noindent

{\it Notation.}

By $S_{n}$, $A_{n}$, $C_{n}$ and $\mu_{n} \simeq C_{n}$, we denote
the symmetric group of
degree $n$,  the alternating group of degree $n$, the cyclc group of order $n$, and the multiplicative group of order
$n$ (in $\mathbf C^{\times}$)
respectively. Then $\mu_{n} = \langle \zeta_{n} \rangle$ where
$\zeta_{n} = \text{exp}\, (2\pi i/n)$. By $C_{n}^{m}$, we denote the direct
product of $m$ copies of $C_{n}$.

By $\text{PGL}(n, q)$ (resp. $\text{PSL}(n, q)$), we denote
the projective linear
group (resp. projective special linear group) of
$\mathbf P(\mathbf F_{q}^{\oplus n})$ over the
finite field $\mathbf F_{q}$ of $q$-elements.

$M_{10}$ is the Mathieu group of degree $10$,
which is defined to be the pointwise stabilizer subgroup of
$\{11, 12\}$ of the
Mathieu group $M_{12}$
of degree $12$ under the natural action of  $M_{12}$ on the twelve
element set
$\{1, 2, \cdots , 12\}$. (See for instance [CS, Chapters 10, 11]).

We write $G = A.B$ when $G$ fits in the exact sequence
$$1 \longrightarrow A \longrightarrow G
\longrightarrow B \longrightarrow 1\,\, .$$
(So, for given $A$ and $B$, there are in general several
isomorphism classes of
groups of the form $A.B$.) For instance $\text{Aut}(A_{6}) = A_{6}.C_{2}^{2}$.
We always regard $A$ as a normal subgroup of $A.B$.

When the exact sequence above splits, $G$ is a semi-direct product
of $A$ and $B$ which we shall write $A:B$. (Again, for given
$A$ and $B$, there are in general several isomorphism classes of
groups of the form $A:B$.)
We always regard $A$ as a normal subgroup of $A:B$ and $B$ as a
subgroup of  $A:B$.

By $[G, G]$ and $\text{Cent}(G)$, we denote the commutator subgroup of
$G$ and the center of $G$ respectively.
\par \vskip 1pc \noindent

{\it Acknowledgement.} We would like
to express our thanks to Professor A. A. Ivanov for his warm encouragement,
and to Professors M. Gizatullin and I. Dolgachev for pointing out
that $A_6$ is also called Valentiner's group and informing us the reference [Va].

\section{The isomorphism classes of $A_{6}.\mu_{4}$}

Let $G$ be a group of the form $A_{6}.\mu_{4}$, i.e. a group
which fits in the exact sequence
$$1 \longrightarrow A_{6} \longrightarrow G
\, \mapright{\alpha}\, \mu_{4} \longrightarrow 1\,\,  .$$

The goal of this section is to determine the isomorphism classes of such $G$.

In order to make our argument clear, we first remark the following:

\begin{lemma}\label{lemma:hom}

\begin{list}{}{
\setlength{\leftmargin}{10pt}
\setlength{\labelwidth}{6pt}
}
\item[(1)] Both $G$ and $\text{Aut}(A_{6})$ have exactly one subgroup which is isomorphic to $A_{6}$.
\item[(2)] $\text{Aut}(A_{6}) \simeq A_{6}.C_{2}^{2}$ has exactly one subgroup
which is isomorphic to $S_{6}$, $\text{PGL}(2, 9)$, and $M_{10}$
respectively. These three groups are mutually non-isomorphic
even as abstract groups. In terms of the
natural (conjugacy) action on $A_{6}$, these three subgroups are also
distinguished as follows: $S_{6}$ switches the
conjugacy classes $5A$ and $5B$ of $A_{6}$, but not $3A$ and $3B$;
$\text{PGL}(2, 9)$ switches the
conjugacy classes $3A$ and $3B$ of $A_{6}$ but not $5A$ and $5B$; and
$M_{10}$ switches both.
(The notation here is from the Atlas; see also the
table in Proposition (3.2).)

\end{list}
So, for instance one can speak of {\it the} subgroup $M_{10}$ of $\text{Aut}(A_{6})$
without any ambiguity.
\end{lemma}
\begin{proof} Both $G$ and $\text{Aut}(A_{6})$ have a normal subgroup $H$ isomorphic to
$A_{6}$. This $H$ is the kernel of the natural homomorphism $G \longrightarrow \mu_{4}$
and $\text{Aut}(A_{6}) \longrightarrow C_{2}^{2}$ respectively. Let $K$
be another subgroup of
$G$ or $\text{Aut}(A_{6})$ isomorphic to $A_{6}$.
Since $K \simeq A_{6}$ is a simple non-commutative group, the homomorphisms
$K \longrightarrow \mu_{4}$ and $K \longrightarrow C_{2}^{2}$ are both trivial.
Thus $H > K$ and hence $K = H$ by $\vert K \vert = \vert H \vert$.

Let $L$ be a subgroup of $\text{Aut}(A_{6})$ being isomorphic to
either $S_{6}$, $\text{PGL}(2, 9)$ or $M_{10}$. In each case $L$
contains a subgroup isomorphic to $A_{6}$. This same $A_{6}$ is
now the unique subgroup of $\text{Aut}(A_{6}) = A_{6}.C_{2}^{2}$
and therefore such $L$ corresponds bijectively to the index
$2$-subgroups of $\text{Aut}(A_{6})/A_{6} = C_{2}^{2}$. Thus,
there are exactly three such $L$. The last two assertions and the
fact that $\text{Aut}(A_{6})$ contains subgroups isomorphic to
$S_{6}$, $\text{PGL}(2, 9)$ and $M_{10}$ are well known. (See for
instance [CS, Chap 10] or Atlas table of $A_{6}$.)

\end{proof}

Since $A_{6}$ is a normal
subgroup of $G$, it follows that $c(g)(a) = gag^{-1} \in A_{6}$ for
$g \in G$ and $a \in A_{6}$.  We have then a natural group homomorphism
$$c : G \longrightarrow\, \text{Aut}\,(A_{6})\,\, ;\,\, g \mapsto c(g)\,\, ;\,\,
c(g)(a) = gag^{-1}\,\, .$$
Set
$$N\,\, :=\,\, c(G)$$
and consider the homomorphism:
$$\tilde{c} := (c, \alpha) : G \longrightarrow N \times \mu_{4}
< \text{Aut}(A_{6})
\times \mu_{4}\,\,
;\,\, g \mapsto (c(g), \alpha(g))\,\, .$$
We define the natural projections:
$$p_{1} : N \times \mu_{4} \longrightarrow N\,\,
;\,\, (x, y) \mapsto x\,\, ,\,\,
p_{2} : N \times \mu_{4} \longrightarrow \mu_{4}\,\,
;\,\, (x, y) \mapsto y\,\, .$$

\begin{lemma}\label{lemma:hom}
\begin{list}{}{
\setlength{\leftmargin}{10pt}
\setlength{\labelwidth}{6pt}
}
\item[(1)] $N$ is either $A_{6}$, $S_{6}$,
$\text{PGL}(2, 9)$ or $M_{10}$.
\item[(2)] $\tilde{c}$ is injective.
\item[(3)] $\alpha^{-1}(\mu_{2}) \simeq A_{6} \times \mu_{2}$
and $\tilde{c}(\alpha^{-1}(\mu_{2})) = A_{6} \times \mu_{2}$ in
$N \times \mu_{4}$. In particular, if $\tilde{h} \in G$ and $c(\tilde{h}) \in
N \setminus A_{6}$, then $\alpha(\tilde{h}) = \pm \zeta_{4}$.
\end{list}
\end{lemma}

\begin{proof} Since $A_{6}$ is simple and non-commutative,
the restriction $c \vert A_{6}$ is injective.
Thus
$$A_{6} = c(A_{6}) < N < \text{Aut}(A_{6}) = A_{6}.C_{2}^{2}\,\, .$$
Moreover, we have that $N \not= A_{6}.C_{2}^{2}$; otherwise
$c : A_{6}.\mu_{4} \simeq A_{6}.C_{2}^{2}$ is an isomorphism and
$\mu_{4} \simeq C_{2}^{2}$ by Lemma (2.1),
a contradiction. Thus,  we obtain (1).

Let $g \in \text{Ker}(\tilde{c})$. Then, $\alpha(g) = 1$, i.e. $g \in A_{6}$. Since $A_{6}$ is simple
and non-commutative, it follows that $\text{Cent}(A_{6}) = \{1\}$. Thus
$g = 1$ by $c(g) = 1$,
and we obtain (2).

Since $\mu_{2}$ is a normal subgroup of $\mu_{4}$ of index $2$ and $\alpha$
is surjective,
it follows that $\alpha^{-1}(\mu_{2})$ is a normal subgroup of $G$ of index $2$. Thus
$[\alpha^{-1}(\mu_{2}) : A_{6}] = 2$. Take $g \in G$ such that
$\alpha(g) = \zeta_{4}$.
Then $\alpha^{-1}(\mu_{2}) = \langle A_{6}, g^{2} \rangle$. Since
$[c(G) : A_{6}] \leq 2$ by (1), we have $c(g^{2}) \in A_{6}$. Thus, $c(g^{2}h^{-1}) = 1$ for some
$h \in A_{6}$. Set $f := g^{2}h^{-1}$. Then $\tilde{c}(f) = (1, -1)$,
whence $\text{ord}\, f = 2$ by (2), and $fa = af$ for each $a \in A_{6}$.
Thus
$$\alpha^{-1}(\mu_{2}) = \langle A_{6}, f \rangle =
A_{6} \times \langle f \rangle \simeq A_{6} \times \mu_{2}\,\, .$$
This implies
$$A_{6} \times \mu_{2} < \tilde{c}(\alpha^{-1}(\mu_{2}))
< N \times \mu_{2}$$
in $N \times \mu_{4}$, and hence $A_{6} \times \mu_{2}
= \tilde{c}(\alpha^{-1}(\mu_{2}))$, because the orders are the same.
\end{proof}

The four candidates for $N$ in Lemma (2.2)(1) give four different
group structures on $G = A_{6}.\mu_{4}$, which are all semi-direct
products, indeed:

\begin{theorem}\label{theorem:abstract group}
There are exactly four possible group structures of $G$, up to
isomorphism. More explicitly, $G$ is isomorphic to the following
subgroup of $N \times \mu_{4}$ corresponding to each of the four
candidates for $N < \text{Aut}(A_{6})$ as in (2.2)(1):

\begin{list}{}{
\setlength{\leftmargin}{10pt}
\setlength{\labelwidth}{6pt}
}
\item[(1)] If $N = A_{6}$, then $G = A_{6} \times \mu_{4}$.
In this case, we set
$\tilde{g} = (1, \zeta_{4})$. Then $\tilde{g}$ is an order $4$ element and
$G = A_{6} \times \langle \tilde{g} \rangle$.
\item[(2)] If $N = S_{6}$, then $G = A_{6}.\langle \tilde{g} \rangle =
A_{6} : \langle \tilde{g} \rangle$, where $\tilde{g} = (g, \zeta_{4})$
and $g = (1, 2) \in S_{6}$.
\item[(3)] If $N = \text{PGL}(2,9)$, then
$G = A_{6}.\langle \tilde{g} \rangle
= A_{6} : \langle \tilde{g} \rangle$, where
$\tilde{g} = (h^{5}, \zeta_{4})$ for some element $h \in \text{PGL}(2, 9)$
of order $10$.
\item[(4)] If $N = M_{10}$, then $G = A_{6}.\langle \tilde{g} \rangle =
A_{6} : \langle \tilde{g} \rangle$, where $\tilde{g} = (g, \zeta_{4})$ and
$g$ is an order $4$-element of $M_{10} \setminus A_{6}$.
\end{list}
In each case, the semi-direct product structure is the natural
one. We denote the groups in (1), (2), (3), (4) by $A_{6}(4)$,
$S_{6}(2)$, $\text{PGL}(2,9)(2)$, $M_{10}(2)$ respectively.
\end{theorem}

\begin{proof} It is clear that the four groups $A_{6}(4)$, $S_{6}(2)$,
$\text{PGL}(2,9)(2)$, $M_{10}(2)$ satisfy the
required conditions,
i.e. $G$ is of the form $A_{6}.\mu_{4}$ and $c(G) = N$.

Let us next show that $G$ is isomorphic to the group $A_{6}(4)$,
$S_{6}(2)$, $\text{PGL}(2,9)(2)$, $M_{10}(2)$, if $N = A_{6}$,
$S_{6}$, $\text{PGL}(2,9)$, $M_{10}$ respectively.

If $c(G) = A_{6}$, then $\tilde{c} : G \longrightarrow A_{6} \times \mu_{4}$ is an isomorphism.
Indeed, $\tilde{c}$ is injective by Lemma (2.2)(2) and
$\vert G \vert = \vert A_{6} \times \mu_{4}
\vert$.

Consider the case $N = S_{6}$ or $\text{PGL}(2, 9)$. Since $c(G) = N$, there is an element
$\tilde{g}$ of $G$ such that $c(\tilde{g}) \in N \setminus A_{6}$,
$\text{ord}(c(\tilde{g})) = 2$ and $\alpha(\tilde{g}) = \pm \zeta_{4}$.
Indeed, by Lemma (2.2)(3), one can take a preimage of $g$ in (2) and (3)
as $\tilde{g}$. Moreover,
$\zeta_{4} \mapsto -\zeta_{4}$ gives a group
automorphism of $\mu_{4}$, the isomorphism class of $G$ does not
depends on the choice of the sign of $\alpha(\tilde{g})$.
Therefore we may also adjust as
$\alpha(\tilde{g}) = \zeta_{4}$ for a chosen $g$.
Since $\tilde{c}$ is injective,
this also implies $\text{ord}(\tilde{g}) = 4$ and consequently
$G = A_{6} : \langle \tilde{g} \rangle$ as claimed.

Finally consider the case $N = M_{10}$. Note that $M_{10} \setminus A_{6}$
has no order $2$
element and order $4$ elements of
$M_{10} \setminus A_{6}$ form one conjugacy class of $M_{10}$.
(See for instance, [Atlas table of $A_{6}$].) Let
$g \in M_{10} \setminus A_{6}$ be an
order $4$ element
and $\tilde{g} \in G$ be an element such that $c(\tilde{g}) = g$.
Then $\alpha(\tilde{g}) = \pm \zeta_{4}$ by Lemma (2.2)(3).
Now, as in
the cases (2) and (3), we may adjust as
$\alpha(\tilde{g}) = \zeta_{4}$ for a chosen $g$, and
$G = A_{6} : \langle \tilde{g} \rangle$ as claimed in (4).

Since $A_{6}$ is the unique subgroup of $G$, the image $N$
of the homomorphism $c$ is uniquely determined by $G$. Thus, the
four groups $A_{6}(4)$, $S_{6}(2)$,
$\text{PGL}(2,9)(2)$, $M_{10}(2)$ are not isomorphic to one another.

Now we are done.

\end{proof}

\section{Some lattice representations}

In this section, we recall some known facts about K3 surface from
[BPV] and about K3 surfaces admitting an $A_{6}$-action from
[KOZ]. For details, please refer to these references and
references therein.

By a $K3$ surface, we mean a
simply-connected compact complex surface $X$ admitting a nowhere
vanishing global holomorphic $2$-form $\omega_{X}$. K3 surfaces form
a $20$-dimensional analytic family. The second
cohomology group $H^{2}(X, \mathbf Z)$ together with its cup product
becomes an even unimodular lattice of index $(3, 19)$ and is
isomorphic to the so-called $K3$ lattice
$$L:= U^{\oplus 3} \oplus
E_{8}^{\oplus 2}\,\, ,$$ where $E_8$ is the negative definite even
unimodular lattice of rank 8. We denote by $S(X)$ the
N\'eron-Severi lattice of $X$. This is a primitive sublattice of
$H^{2}(X, \mathbf Z)$ generated by the (first Chern) classes of
line bundles. The rank of $S(X)$ is called the Picard number of
$X$ and is denoted by $r(X)$. We have $0 \leq r(X) \leq 20$.
Denote by $T(X)$ the transcendental lattice of $X$, i.e. the
minimum primitive sublattice whose $\mathbf C$-linear extension
contains the class $\omega_{X}$, or equivalently $T(X) =
S(X)^{\perp}$ in $H^{2}(X, \mathbf Z)$. If $X$ is projective, then
$S(X)$ is of index $(1, r(X)-1)$ (and vice versa), $S(X) \cap T(X)
= \{0\}$ and $S(X) \oplus T(X)$ is a finite-index sublattice of
$H^{2}(X, \mathbf Z)$.

Let $(X, G, \rho)$ be a triplet consisting of a $K3$ surface $X$,
a finite group $G$ and a faithful action $\rho : G \times X
\longrightarrow X$. Then $G$ has a $1$-dimensional representation
on $H^{0}(X, \Omega_{X}^{2}) = \mathbf C \omega_{X}$ given by
$g^{*}\omega_{X} = \alpha(g)\omega_{X}$, and we have an exact
sequence, called the basic sequence:
$$1 \longrightarrow G_{N} := \text{Ker}\, \alpha \longrightarrow G
\mapright{\alpha}
\mu_{I} \longrightarrow 1\,\,  .$$

The importance of the basic sequence was first noticed by Nikulin [Ni].
We call $G_{N}$ the symplectic
part and $\mu_{I} := \langle \zeta_{I} \rangle$ (resp. $I$),
the transcendental
part (resp. the transcendental value) of the action $\rho$. We
note that if $A_{6}.\mu_{4}$ acts faithfully on a $K3$ surface
then $G_{N} \simeq A_{6}$ and the transcendental part is
isomorphic to $\mu_{4}$. This follows from the fact that $A_{6}$
is simple and also maximum among all finite groups acting on a
$K3$ surface symplectically. This is a result of Mukai [Mu]. We also note
that
$X$ is projective if $I \geq 2$ [Ni].

We say that 2 triplets $(X, G, \rho)$ and $(X', G', \rho')$ are isomorphic if
there are a group isomorphism $f : G' \simeq G$ and an isomorphism
$\varphi : X' \simeq X$ such that the following diagram commutes:
$$
\begin{matrix}
G \times X & \mapright{\rho} & X \cr
\mapupleft{f \times \varphi \; } &  & \mapup{\varphi} \cr
G' \times X' & \mapright{\rho'} & X' \cr
\end{matrix}
$$

The main result of [KOZ] is the following:
\begin{theorem} \label{theorem:prep}
Let $G$ be a finite group acting faithfully on a $K3$ surface $X$. Assume that
$A_{6} < G$ and $I \geq 2$, where $I$ is the transcendental value. Then we have
$G_{N} = A_{6}$, $I = 2$ or $4$, and ${\rm rank}\, L^{G_{N}} = 3$.
In particular,
$S(X)^{G_{N}} =
\mathbf Z H$, where $H$ is an ample class, and ${\rm rank}\, T(X)
= 2$.

Consider the case having maximum $I = 4$, i.e. the case where $G$,
which we shall denote by $\tilde{A}_{6}$, is a group of the form
$A_{6}.\mu_{4}$. Then there is a unique triplet $(F,
\tilde{A}_{6}, \rho_{F})$ consisting of a $K3$ surface $F$ and a
faithful group action $\rho_{F} : \tilde{A}_{6} \times F
\longrightarrow F$ of $\tilde{A}_{6}$ on $F$ up to isomorphism. In
particular, the isomorphism class of $\tilde{A}_{6}$ is unique.
Moreover the K3 surface $F$ has Picard number $20$ and is uniquely
characterized by the following equivalent conditions:
\begin{list}{}{
\setlength{\leftmargin}{10pt}
\setlength{\labelwidth}{6pt}
}
\item[(1)]
$F$ is the $K3$ surface whose transcendental lattice
$T(F) = \mathbf Z \langle t_{1}, t_{2} \rangle$ has the intersection matrix
$$\begin{pmatrix}
6&0\cr
0&6
\end{pmatrix}\,\, .$$
\item[(2)] $F$ is the minimal resolution of
the double cover $\overline{F}$ of the (rational) elliptic modular surface $E$
with level $3$ structure. The double cover is branched along two
of a total of $4$ singular fibres of the same type $I_3$ and $\overline{F}$ has
$6$ ordinary double points.
\item[(3)]
$F$ is the minimal resolution of the surface in $\mathbf P^{1} \times
\mathbf P^{2}$ given by the following equation, where
$([S:T],[X:Y:Z])$ are coordinates of $\mathbf P^{1} \times \mathbf
P^{2}$:
$$S^2(X^3+Y^3+Z^3) - 3(S^2+T^2)XYZ = 0.$$
\end{list}
\end{theorem}

In the course of proof, we have also shown the following fact.
This will also be crucial in the next section and will be proved
again here:

\begin{proposition} \label{proposition:irred} Let $X$ be a K3 surface of Picard number
$20$. Assume that $X$ admits a faithful action of $G_{N} = A_{6}$. Then,
as $A_{6}$-modules, one
has the irreducible decomposition
$$S(X) \otimes \mathbf C = \chi_{1} \oplus \chi_{2} \oplus \chi_{3}
\oplus \chi_{6}\,\, .$$
In this description, we used
Atlas notation for irreducible characters/representations of $A_{6}$ as in
the Table below.
\end{proposition}

\par \vskip 3pc


\begin{tabular}{|l|l|l|l|l|l|l|l|l|r|} \hline \hline
              &\vline \, 1A  &\vline \, 2A  &\vline \, 3A  &\vline \, 3B  & \vline \, 4A
&\vline \, 5A & \vline \, 5B
\\ \hline \hline

$\chi_1$ \, &\vline \, 1  &\vline \, 1 &\vline \, 1 &\vline \, 1  & \vline \, 1
&\vline \, 1 &\vline \, 1
\\ \hline \hline

$\chi_2$ \, &\vline \, 5  &\vline \, 1 &\vline \, 2   &\vline \, -1  & \vline \, -1
&\vline \, 0  &\vline \, 0
\\ \hline \hline

$\chi_3$  &\vline \, 5    &\vline \, 1  & \vline \, -1  &\vline \, 2 &\vline \, -1
&\vline \,  0  &\vline \, 0
\\ \hline \hline

$\chi_4$  &\vline \, 8    &\vline \, 0  &\vline \, -1   &\vline \,
-1 &\vline \, 0 &\vline \,  $({1-\sqrt{5}})/2$  &\vline \,
$({1+\sqrt{5}})/{2}$
\\ \hline \hline

$\chi_5$  &\vline \, 8    &\vline \, 0  &\vline \, -1 &\vline \,
-1 & \vline \, 0 &\vline \,  $({1+\sqrt{5}})/{2}$  &\vline \,
$({1-\sqrt{5}})/{2}$
\\ \hline \hline

$\chi_6$  &\vline \, 9    &\vline \ 1   &\vline \, 0  &\vline \, 0
& \vline \, 1 &\vline \, -1  &\vline \, -1
\\ \hline \hline

$\chi_7$  &\vline \, 10   &\vline \ -2  &\vline \, 1  &\vline \, 1  &\vline \, 0
&\vline \,  0  &\vline \, 0
\\ \hline \hline

\end{tabular}

\begin{proof} Since many ingredients and some idea in the next section
have already appeared in the proof
of this proposition, we prove this proposition for reader's convenience.

Recall that the order structure of $A_{6}$ is as follows:
$$
\begin{tabular}{*{8}{|c}{|}}
\hline
\text{order[conjugacy class]} & 1 [1A] & 2 [2A] & 3 [3A] & 3 [3B] & 4 [4A]
& 5 [5A] & 5 [5B]\\
\hline
\text{cardinality}  & 1 & 45 & 40 & 40 & 90 & 72 & 72\\
\hline
\end{tabular}
$$
Moreover, by [Ni], the number of the fixed points of the symplectic
action is as follows.
$$
\begin{tabular}{*{9}{|c}{|}}
\hline
$ord(g)$ & 1 & 2 & 3 & 4 & 5 & 6 & 7 & 8 \\
\hline
$\vert X^{g} \vert$ & $X$ & 8 & 6 & 4 & 4 & 2 & 3 & 2 \\
\hline
\end{tabular}
$$
Set $\tilde{H}(X, \mathbf Z) = H^{0}(X, \mathbf Z) \oplus
H^{2}(X, \mathbf Z) \oplus H^{4}(X, \mathbf Z)$.
Now, by applying the topological Lefschetz fixed point formula for
$G_{N} = A_{6}$,
we calculate that
$$\text{rank}\, \tilde{H}(X, \mathbf Z)^{A_{6}}  =
\frac{1}{\vert A_{6} \vert}
\sum_{g \in A_{6}} \text{tr}(g^{*} \vert \tilde{H}(X, \mathbf Z))$$
$$ =
\frac{1}{360}(24 + 8 \cdot 45 + 6 \cdot 80 + 4 \cdot 90 + 4 \cdot 144 )
= 5\,\, .$$
Since $G_{N}$ is trivial on $H^{0}(X, \mathbf Z)$, $H^{4}(X, \mathbf Z)$ and
$T(X)$,
one has $S(X)^{A_{6}} = \mathbf Z H$ and $H$ is an ample primitive
class of $X$.
(Here we used the fact that any K3 surface of Picard number $20$ is projective, because
$S(X) \otimes \mathbf R = H^{1,1}(X, \mathbf R)$ and therefore $S(X)$ is of signature $(1, 19)$
and the fact that if $h$ is an ample class of $X$ and $G$ is a finite group
acting on $X$,
then $\sum_{g \in G}g^{*}h$ is also an ample class which is invariant under $G$.)

Thus the irreducible decomposition of $S(X)$ by $A_{6}$ must be of the
following form:

$$S(X) \otimes \mathbf C = \chi_{1} \oplus a_{2}\chi_{2} \oplus a_{3}\chi_{3} \oplus
a_{4}\chi_{4} \oplus a_{5}\chi_{5} \oplus a_{6}\chi_{6} \oplus
a_{7}\chi_{7}\,\, ,$$ where $a_{i}$ are non-negative integers. Let
us determine $a_{i}$'s. As in (2), using the topological Lefschetz
fixed point formula and the fact that $\text{rank}\, T(X) = 2$, we
have
$$\chi_{\text{top}}(X^{g}) = 4 + \text{tr}(g^{*} \vert S(X))$$
for $g \in A_{6}$. Running $g$ through the $7$-conjugacy classes of $A_{6}$
and calculating both sides based on Nikulin's table and the character table
above, we obtain the following system of equations:

$$20 = 1 + 5(a_{2} + a_{3}) + 8(a_{4} + a_{5}) + 9a_{6} + 10a_{7}\,\, ,$$
$$4 = 1 +  (a_{2} + a_{3}) + a_{6} - 2a_{7}\,\, ,$$
$$2 = 1 +  (2a_{2} - a_{3}) - (a_{4} + a_{5}) + a_{7}\,\, ,\,\,
2 = 1 +  (-a_{2} +2a_{3}) - (a_{4} + a_{5}) + a_{7}\,\,  ,$$
$$0 = 1 -  (a_{2} + a_{3}) + a_{6}\,\, ,$$
$$0 = 1 +  (\frac{1 - \sqrt{5}}{2}a_{4} + \frac{1 + \sqrt{5}}{2} a_{5}) - a_{6}\,\, ,\,\,
0 = 1 + (\frac{1 + \sqrt{5}}{2}a_{4} + \frac{1 - \sqrt{5}}{2} a_{5})  - a_{6}\,\,  .$$

Now, we get the result by solving this system of Diophantine equations.
\end{proof}

\section{Determination of the isomorphism class of $\tilde{A}_{6}$}

Let $\tilde{A}_{6}$ be a group of the form $A_{6}.\mu_{4}$ which
can act on a K3 surface $X$. Among four candidates $A_{6}(4)$,
$S_{6}(2)$, $\text{PGL}(2,9)(2)$, and $M_{10}(2)$ for
$\tilde{A}_{6}$ (Theorem (2.3)), only one is isomorphic to
$\tilde{A}_{6}$ (Theorem (3.1)).

The aim of this section is to determine the isomorphism class of
this $\tilde{A}_{6}$:

\begin{theorem} \label{theorem:gp} $\tilde{A}_{6}$ is isomorphic
to the group $M_{10}(2)$.
\end{theorem}

\begin{proof} Set $G = \tilde{A}_{6}$. It suffices to show that $G$ is not isomorphic to $A_{6}(4)$, $S_{6}(2)$,
$\text{PGL}(2,9)(2)$. Suppose to the contrary that $G$ is
isomorphic to one of these three groups. Let $\tilde{g} \in G$ be
an order $4$ element chosen in Theorem (2.3). Set $\iota :=
\tilde{g}^{2}$. Then, by the description of $\tilde{g}$, we have
$G = A_{6} : \langle \tilde{g} \rangle$ and $c(\iota) = 1$, i.e.
$\iota a = a \iota$ for all $a \in A_{6}$, and $\alpha(\iota) =
-1$. Thus $\iota^{*}\omega_{X} = -\omega_{X}$ on a target K3
surface $X$.

First we prove the result below to get some geometric constraint of
the pair $(X, \tilde{A}_{6})$.
\begin{lemma} \label{lemma:geom}

\begin{list}{}{
\setlength{\leftmargin}{10pt}
\setlength{\labelwidth}{6pt}
}
\item[(1)]
$\chi_{\text{top}}(X^{\iota}) \leq 0$. Moreover,
$\chi_{\text{top}}(X^{\iota}) = 0$
if and only if $X^{\iota} = \emptyset$.
\item[(2)] Let $\sigma$ be an element of $A_{6}$. Assume that the order of
$\sigma$
is either $3$ or $5$. Then $\chi_{\text{top}}(X^{\iota \sigma}) \geq 0$.
\end{list}
Here $\chi_{\text{top}}(X^{a})$ is the topological Euler number of the fixed locus $X^{a}$
of $a \in G$.
\end{lemma}

\begin{proof} Since $\iota^{*}\omega_{X} = -\omega_{X}$, the action of
$\iota$ on $X$
is locally linearlizable
at $P \in X^{\iota}$ as $$\begin{pmatrix}
1&0\cr
0&-1
\end{pmatrix}\,\, .$$

Set $C := X^{\iota}$.  If it is empty, then $\chi_{\text{top}}(C)
= 0$. Assume that $C \not= \emptyset$. Then $C$ is a smooth curve,
possibly reducible. Since $\iota a = a \iota$ for each $a \in
A_{6}$ (by $c(\iota) = 1$), the curve $C$ is stable under the
action of $A_{6}$. Thus the class $[C] \in S(F)$ is
$A_{6}$-invariant. Since $S(X)^{A_{6}} = \mathbf Z H$ and $H$ is
an ample primitive class, $C$ is also an ample class. Therefore,
$C$ is connected and hence irreducible by the smoothness. Since
$(C^{2}) > 0$,  we have by the adjunction formula
$$\chi_{\text{top}}(C) = 2 - 2g(C) = - (K_{X} + C. C) = -(C^{2}) < 0\,\, .$$

Let us show (2). Assume that $\sigma$ is of order $p$, where $p$ is either $3$ or $5$.
Set $\tau = \iota \sigma$. Then $X^{\tau} \subset X^{\tau^{2}}$. We have
$\tau^{2} = \sigma^{2} \in A_{6} = G_{N}$ (and is of order $p$) by
$\iota \sigma = \sigma \iota$.
Thus, the set $X^{\tau^{2}}$ is a finite set as in Nikulin's table in
Proposition (3.2). Hence $X^{\tau}$ is also a finite set, and therefore
$\chi_{\text{top}}(X^{\tau}) \geq 0$.

\end{proof}

Let us recall the irreducible decomposition of $S(X)$ in
Proposition (3.2). Since $X$ has an ample $G$-invariant class, we
have $\tilde{g}^{*}(H) = H$ and hence $\tilde{g}^{*} \vert
\chi_{1} = id$. Since $G = A_{6} : \langle \tilde{g} \rangle$, we
have also $\tilde{g}^{*}(\chi_{6}) = \chi_{6}$ and either
$\tilde{g}^{*}(\chi_{2}) = \chi_{2}$ and $\tilde{g}^{*}(\chi_{3})
= \chi_{3}$ or $\tilde{g}^{*}(\chi_{2}) = \chi_{3}$ and
$\tilde{g}^{*}(\chi_{3}) = \chi_{2}$. Since $\iota =
\tilde{g}^{2}$, we have $\iota^{*}(\chi_{i}) = \chi_{i}$ for each
$i = 1$, $2$, $3$, $6$. Since $\iota a = a \iota$ for all $a \in
A_{6}$, it follows that $\iota^{*} \vert \chi_{i}$ are scalar
multiplications by Schur's lemma. Moreover, since $\iota$ is of
order $2$, we have
$$\iota^{*} \vert \chi_{1} = 1\,\, , \,\,\iota^{*}  \vert \chi_{i} =
(-1)^{n_{i}}id_{\chi_{i}}\,\,$$
for some $n_{i} \in \mathbf Z$ for each $i = 2$, $3$, $6$. We have also
that $\iota^{*} \vert H^{0}(X) \oplus H^{4}(X) = id$,
$\iota^{*} \vert T(X) = -id_{T(X)}$ and $\text{rank}\, T(X) = 2$.
Thus, by the topological Lefschetz
formula, we obtain
$$\chi_{\text{top}}(X^{\iota}) = 1 +
5 \cdot ((-1)^{n_{2}} + (-1)^{n_{3}}) +
9 \cdot (-1)^{n_{6}}\,\, ---\,\, (*)\,\, .$$

The value $\chi_{\text{top}}(X^{\iota})$ must be non-positive by
Lemma (4.2)(1), whence
$$((-1)^{n_{2}}, (-1)^{n_{3}}, (-1)^{n_{6}})$$
must be
either one of:

$$(-1, 1, -1)\,\, ,\,\, (1, -1, -1)\,\, ,\,\, (-1, -1, -1)\,\, ,\,\,
(-1, -1, 1)\,\, .$$

Consider first the case $(-1, 1, -1)$ (resp. $(1, -1, -1)$).
Take an order $3$ element
$\sigma$ of $A_{6}$ from the conjugacy class 3A (resp. 3B). Note that
$(\iota \sigma)^{*} \vert H^{0}(X) \oplus H^{4}(X) = id$,
$(\iota \sigma)^{*} \vert T(X) = -1$ and $\text{rank}\, T(X) = 2$. Then by the topological Lefschetz
formula and by the character table, we calculate
$$\chi_{\text{top}}(X^{\iota \sigma}) =
\text{tr} ((\iota \sigma)^{*} \vert S(X)) = 1 - 2 -1 + 0 = -2 < 0\,\, ,$$
a contradiction to Lemma (4.2)(2).

Consider next the case $(-1, -1, -1)$. Since
$\tilde{g}^{*}(\chi_{6}) = \chi_{6}$ and $\iota^{*} \vert \chi_{6}
= -id_{\chi_{6}}$, it follows that
the eigenvalues of $\tilde{g}^{*} \vert \chi_{6} = \pm \zeta_{4}$
and $\text{tr}(\tilde{g}^{*} \vert \chi_{6}) = (9-2n)\zeta_{4}$,
where $n$ is the multiplicity
of $-\zeta_{4}$. Note that
$\text{tr}(\tilde{g}^{*} \vert T(X)) =  \zeta_{4} -\zeta_{4} = 0$.

So, if $\tilde{g}^{*}(\chi_{2}) = \chi_{3}$, then
$\text{tr}(\tilde{g}^{*} \vert \chi_{2} \oplus \chi_{3}) = 0$. Thus
$$\chi_{\text{top}}(X^{\tilde{g}}) = 2 + \text{tr}(\tilde{g}^{*} \vert S(X))
= 3 + (9-2n)\cdot \zeta_{4} \not\in \mathbf Z\,\, ,$$
a contradiction to the obvious fact that
$\chi_{\text{top}}(X^{\tilde{g}}) \in \mathbf Z$.

If $\tilde{g}^{*}(\chi_{2}) = \chi_{2}$, then for the same reason
as above, the
eigenvalues of
$$\tilde{g}^{*} \vert \chi_{2} \oplus \chi_{3} \oplus \chi_{6}$$
are $\pm \zeta_{4}$. Let $n$ be the multiplicity
of $-\zeta_{4}$. Then, since
$$\text{dim}\, \chi_{2} \oplus \chi_{3} \oplus \chi_{6}= 19\,\, ,$$
we have
$$\chi_{\text{top}}(X^{\tilde{g}}) =
2 + \text{tr}(\tilde{g}^{*} \vert S(X)) = 3 + (19-2n)\zeta_{4}
\not\in \mathbf Z\,\, ,$$
again a contradiction.

Finally consider the case $(-1, -1, 1)$.

Let us first treat the cases where $G$ is isomorphic to $A_{6}(4)$ or
$S_{6}(2)$. By the formula (*), we have $\chi_{\text{top}}(X^{\iota}) = 0$
and therefore $X^{\iota} = \emptyset$ by Lemma (4.2)(1).
Let $\tau := (456)$ in $A_{6}$.
Then, by the shape of $\tilde{g}$, we have $\tau \tilde{g} = \tilde{g} \tau$
in $G$. Then $\tilde{g}$ acts on $X^{\tau}$. Since
$X^{\tau}$ is a $6$ element set (see Nikulin's table in Section 3)
and $\tilde{g}$ is of order $4$, it follows that $\iota (= \tilde{g}^{2})$
fixes at least two points in $X^{\tau}$. In particular,
$X^{\iota} \not= \emptyset$, a contradiction.

Let us next consider the case where $G$ is isomorphic to $\text{PGL}(2,9)(2)$.

Set $V := (S(X) \otimes \mathbf Q)^{\iota^{*}}$. Then $V$ is $A_{6}$-stable by
$\iota a = a \iota$ for all $a \in A_{6}$. By
$(n_{1}, n_{2}, n_{3}) = (-1, -1, 1)$
and by the fact that the action $\iota^{*}$ is defined over
$S(X) \otimes \mathbf Q$,
we then have
$$V \otimes \mathbf C = \chi_{1} \oplus \chi_{6}\,\, .$$
Recall that $\chi_{1} = \mathbf C H$, where $H$ is an ample
$A_{6}$-invariant class. Consider the orthogonal complement $L :=
H_{V}^{\perp}$ of $H$ in $V$. Then $L$ is also $A_{6}$-stable and
satisfies
$$L \subset S(X) \otimes \mathbf Q\,\, \text{and}\,\,
L \otimes \mathbf C = \chi_{6}\,\, .$$
Since $c(\tilde{g}) = g$ switches the conjugacy classes $3A$ and $3B$ of
$A_{6}$, while $\chi_{2}(3A) = 2 \not= -1 = \chi_{2}(3B)$, we have
$$\tilde{g}^{*}(\chi_{2}) = \chi_{3}\,\, ,\,\,
\tilde{g}^{*}(\chi_{3}) = \chi_{2}\,\, ,$$
and therefore
$$\text{tr}\, (\tilde{g}^{*} \vert \chi_{2} \oplus \chi_{3})\, = \, 0\,\, .$$
Moreover, since $\iota = \tilde{g}^{2}$ and $\iota^{*} \vert \chi_{6} = id$,
we have a matrix representation
$$\tilde{g}^{*} \vert \chi_{6} = I_{9-s} \oplus (-I_{s})$$
under certain {\it rational} basis $\langle u_{i}
\rangle_{i=1}^{9}$ of $L$. This is because $\tilde{g}^{*}$ is
defined over $S(X)$ (so that $\tilde{g}^{*} \vert \chi_{6}$ is
defined over $L$) and the eigenspace decomposition of
$\tilde{g}^{*} \vert \chi_{6}$ of eigenvalues $\pm 1 \in \mathbf
Q$ are also rationally defined. We have also that
$$\text{tr}\,(\tilde{g}^{*} \vert \chi_{6}) = 9 - 2s\,\, .$$
Note also that
$$\text{tr}(g^{*} \vert T(X)) = \zeta_{4} - \zeta_{4} = 0\,\, ,$$
$$\text{tr}(g^{*} \vert H^{0}(X) \oplus H^{4}(X)) = 1 + 1 = 2\,\, .$$
Thus, by the topological Lefschetz fixed point formula, one calculates
$$\chi_{\text{top}}(X^{\tilde{g}}) = 2 + 1 + (9 - 2s) = 12 - 2s\,\, .$$
Since $X^{\iota} = \emptyset$, we have $X^{\tilde{g}} = \emptyset$
as well by $\iota = \tilde{g}^{2}$. Therefore,
$\chi_{\text{top}}(X^{\tilde{g}}) = 0$ and $s = 6$.

Let us recall that $\tilde{c}(\tilde{g}) = (h^{5}, \zeta_{4})$ for some
element $h \in \text{PGL}(2,9)$ of order $10$. Then, $h^{2} \in A_{6}$
(by $[\text{PGL}(2,9) : A_{6}] = 2$) and therefore there is an element
$\sigma \in A_{6} (< G)$ such that $\tilde{c}(\sigma) = (h^{2}, 1)$.
Using the injectivity of $\tilde{c}$, we see that
$\sigma$ is an element of order $5$ and satisfies
$\sigma \tilde{g} = \tilde{g} \sigma$. Note that $\sigma^{*}$
is defined over $S(X)$ and therefore $\sigma^{*} \vert \chi_{6}$
is defined over $L$. Thus, under the same {\it rational}
basis
$\langle u_{i} \rangle_{i=1}^{9}$ of $L$, we have a
{\it rational}
matrix representation
$$\sigma^{*} \vert \chi_{6} = A \oplus B\,\, ,$$
where $A \in \text{GL}(3, \mathbf Q)$ and $B \in \text{GL}(6, \mathbf Q)$.

Since $\sigma \in A_{6}$ is of order $5$, we have $\text{tr}(\sigma^{*} \vert \chi_{6}) = -1$ by the definition of $\chi_{6}$. Thus
$$\text{tr} A + \text{tr} B = -1\,\, .$$
On the other hand, since $\sigma$ is of order $5$, the eigenvalues of $A$
and $B$ are all root of $5$. Since $A$ and $B$ are rational matrices
of orders $3$ and $6$ and since $\varphi(5) = 4$, we have then
$$\text{tr}\, A = 1 + 1 + 1 = 3\,\, ,$$
and either
$$\text{tr}\, B = 1 + 1 + \sum_{i = 1}^{4} \zeta_{5}^{i} = 1\,\, ,$$
or
$$\text{tr}\, B = 1 + 1 + 1 + 1 + 1 + 1 = 6\,\, .$$
However, then
$$\text{tr} A + \text{tr} B = 4\,\, \text{or}\,\, 9\,\, ,$$
a contradiction to the previous equality.

This completes the proof of Theorem (4.1).

\end{proof}

\begin{remark}
In the above proof, the last case $(-1, -1, 1)$ can be ruled out
geometrically. Indeed, in this case we have $X^{\iota} =
\emptyset$ and hence $X^{\tilde{g}} = \emptyset$. This means that
the surface $X$ admits a free action of a cyclic group of order 4,
which is impossible because no K3 surface may admit such an action
(the algebraic Euler number of a K3 surface is equal to 2, not
divisible by 4).
\end{remark}


\end{document}